\documentclass[12pt]{amsart}

\usepackage{latexsym}
\usepackage{amsmath}
\usepackage{amssymb}
\usepackage{amscd}
\usepackage{multicol}
\usepackage{multicol}
\usepackage[cmtip,matrix,arrow]{xy}

\newtheorem{theorem}{\bf Theorem}[section]

\newtheorem{proposition}[theorem]{\bf Proposition}
\newtheorem{lemma}[theorem]{\bf Lemma}

\newtheorem{definition-theorem}[theorem]{\bf Theorem-Definition}

\def\Q{\mathbb{Q}}
\def\Z{\mathbb{Z}}
\def\P{\mathbb{P}}
\def\bR{\mathbb{R}}
\def\bC{\mathbb{C}}
\def\H{\mathbb{H}}
\def\S{\mathcal{S}}

 \def\k{\mathfrak{k}}

\def\d{\frak{d}}

\def\I{{\mathcal I}}

\def\V{{\mathfrak{V}}}

\author{Augustin-Liviu Mare}
\address{
Department of Mathematics and Statistics\\ University of
Regina \\ College West 307.14 \\ Regina SK, Canada S4S
0A2}
\email{mareal@math.uregina.ca}
\date{\today}
\title[Quaternionic flag manifolds]{Equivariant cohomology of quaternionic flag manifolds}

\begin{document}
\begin{abstract} The main result of the paper is a Borel type description of the $Sp(1)^n$-equivariant cohomology ring 
of the manifold $Fl_n(\H)$ of all complete flags in $\H^n$.  To prove this, we obtain   a Goresky-Kottwitz-MacPherson type description of that ring.
\end{abstract}
\maketitle

\section{Introduction}

In this paper we study the quaternionic flag manifold $Fl_n(\H)$, which is the space of
all sequences $(V_1,\ldots,V_n)$ where $V_{\nu}$ is a $\nu$-dimensional quaternionic vector subspace (that is, left $\H$-submodule)  of $\H^n$, for $1\le \nu \le n$, such that $V_1\subset V_2\subset \ldots \subset V_n$. We can see that 
$Fl_n(\H)$ has a transitive action of the symplectic group $Sp(n)$, with stabilizer
$$K:=Sp(1)^n.$$  We are  interested in the  $K$-equivariant cohomology\footnote{All cohomology rings in this paper will be with coefficients in $\Z$ (unless otherwise specified).}  of $Fl_n(\H)$.  
We describe this ring in terms of  the  canonical vector bundles ${\mathfrak V}_{\nu}$   over $Fl_n(\H)$, where
$0\le \nu \le n$. More precisely, we prove the following theorem.
\begin{theorem}\label{main} One has the ring isomorphism
$$H^*_K(Fl_n(\H),\Z) \simeq \Z[x_1,\ldots,x_n,u_1,\ldots,u_n]
/\langle (1+x_1)\ldots (1+x_n)=(1+u_1)\ldots (1+u_n)\rangle,$$
where $x_{\nu}=e_K(\V_{\nu}/\V_{\nu-1})\in H_K^4(Fl_n(\H),\Z)$ is the $K$-equivariant Euler class of 
$\V_{\nu}/\V_{\nu-1}$ and $u_1,\ldots,u_n$ are copies of the generator of $H^4(BSp(1),\Z)$, 
so that 
$$H^*(BK,\Z)=H^*(BSp(1)^n,\Z)=\Z[u_1,\ldots,u_n].$$
\end{theorem}

The strategy used  is as follows. First we prove the Kirwan injectivity type result for the $K$ action on $Fl_n(\H)$; that is, the restriction map
$\imath^*: H^*_K(Fl_n(\H)) \to H^*_K(Fl_n(\H)^K)$ is injective, where  $Fl_n(\H)^K$ denotes the fixed point set of the $K$ action. Then we prove the Goresky-Kottwitz-MacPherson (shortly GKM) type characterization of the image of $\imath^*$. Finally we compare what we have obtained with the GKM picture  of the $T$-equivariant cohomology of 
$Fl_n(\bC)$, where $T=(S^1)^n$ is the maximal torus of $U(n)$. We observe  that there exists an abstract isomorphism between the two cohomology rings, which doubles the degrees.   At the end, we use the known Borel-type description of $H^*_T(Fl_n(\bC))$.

\noindent {\bf Remarks.} 1. The  proof of the GKM presentation of $H^*_K(Fl_n(\H))$ 
we will give in section 3 uses the methods of Tolman and Weitsman [To-We] (see also Harada and Holm [Ha-Ho, section 2]).

2. A presentation of the usual cohomology ring of $Fl_n(\H)$ has been described by us in [Ma1, section 3]
(the result was originally proved by Hsiang, Palais, and Terng [Hs-Pa-Te]). 
Theorem 1.1 gives the equivariant ``deformation" of that presentation.   

\noindent {\bf Acknowledgements.} I wanted to thank  Megumi Harada and Tara Holm for reading a previous version of the paper and making some excellent suggestions. 
I also thank Jost Eschenburg for discussions about the topics of the paper.

\section{The quaternionic flag manifold}\label{second}

The goal of this section is to give a few alternative presentations of the manifold 
$Fl_n(\H)$, which will be used later. 

Let  $$\H=\{a+bi+cj+dk \ | \ a,b,c,d \in \bR \}$$ be the skew field of quaternions. The space $\H^n$ is a $\H$ module with respect to the scalar multiplication from the left. We equip it with the scalar product
$(  \  , \  )$ given by 
$$( h,k)=\sum_{\nu=1}^nh_{\nu}\bar{k}_{\nu},
$$ for all $h=(h_1,\ldots,h_n)$, $k=(k_1,\ldots,k_n)$ in $\H^n$. 
Any linear transformation of $\H^n$ is described by a matrix 
$A\in {\rm Mat}^{n\times n} (\H)$ according to the formula
$$Ah:=h\cdot A^*,$$
where $h=(h_1,\ldots,h_n)\in \H^n$. Here $\cdot$ denotes the matrix multiplication and the superscript $*$ indicates the transposed conjugate of a matrix. We denote by $Sp(n)$ the group of linear transformations $A$ of $\H^n$ with the property that
$$( A.h,A.k) =( h,k),$$
for all $h,k\in \H^n$. Alternatively, $Sp(n)$ consists of all $n\times n$ matrices $A$
with entries in $\H$ with the property that
$A\cdot A^*=I_n$.

The flag manifold $Fl_n(\H)$ defined in the introduction can also be described as the space of all sequences
$(L_1,\ldots,L_n)$ of 1-dimensional $\H$-submodules of $\H^n$ such that 
$L_{\nu}$ is perpendicular to $L_{\mu}$ for all $\mu,\nu\in \{1,2,\ldots,n\}$,
$\mu\neq \nu$. The group $Sp(n)$ acts transitively on the latter space.     
Indeed, let us consider the canonical basis $e_1=(1,0,\ldots,0),\ldots,
e_n=(0,\ldots,0,1)$ of $\H^n$; if $h^1,\ldots,h^n$ is an orthonormal system in $\H^n$, then the matrix $A$ whose columns are $(h^1)^*,\ldots,(h^n)^*$ is in $Sp(n)$ and satisfies $Ae_{\nu}=h^{\nu}$, $1\le \nu \le n$. The $Sp(n)$ stabilizer of the flag $(\H e_1,\ldots,\H e_n)$ consists of diagonal matrices, that is, it is equal to $Sp(1)^n$. In this way we obtain the identification 
\begin{equation}\label{spn} Fl_n(\H)=Sp(n)/Sp(1)^n.\end{equation} 

Yet another presentation of the quaternionic flag manifold can be obtained by   
considering the conjugation action of $Sp(n)$ on the space $${\mathcal H}_n:=\{X\in {\rm Mat}^{n\times n}(\H) \ | \ X=X^*, \ {\rm Trace}(X)=0\}  .$$
 We pick $n$ distinct real numbers
 $r_1,\ldots,r_n$ with $r_1<r_2<\ldots <r_n$ and $\sum_{\nu =1}^nr_{\nu}=0$, and consider the orbit of the diagonal matrix ${\rm Diag}(r_1,\ldots,r_n)$.  For any element $X$ of the orbit there exist mutually orthogonal lines  $L_1,\ldots,L_n$ such that $X|_{L_{\nu}}$ is the multiplication by $r_{\nu}$,
 for all $1\le \nu \le n$. This gives the identification
 \begin{equation}\label{flnh} Fl_n(\H)=Sp(n).{\rm Diag} (r_1,\ldots,r_n).\end{equation}

Finally, we also mention that $Fl_n(\H)$ is an $s$-orbit or a real flag manifold.
More precisely, it is a principal isotropy orbit of the symmetric space $SU(2n)/Sp(n)$. It turns out that the isotropy representation of the latter space is just the conjugation action of $Sp(n)$ on ${\mathcal H}_n$ mentioned above. Via this identification, the metric on
${\mathcal H}_n$ turns out to be the standard one, given by
$$\langle X,Y\rangle :={\rm Re}({\rm Trace}(XY)),$$
$X,Y\in {\mathcal H}_n$. 
 Moreover, a maximal abelian subspace of ${\mathcal H}_n$ is the space $\d$ of all diagonal matrices with real entries and trace 0.  For the details, we refer the reader to  [Ma2, Example 5.4].
Let us consider an element $A={\rm Diag}(a_1,\ldots,a_n)$ of $\d$,
 where  $a_{\nu}$ are real numbers such that $a_1<a_2<\ldots <a_n$ and $\sum_{\nu=1}^na_{\nu}=0$. 
The corresponding height function 
$$h_A(X):=\langle A,X\rangle,$$ $X\in Fl_n(\H)$, will be an important instrument in our paper. 
We will use the following result, a proof of which can be found in the last section of our paper.
\begin{proposition}\label{morsecell} \begin{itemize}
\item[(i)] The critical set 
${\rm Crit}(h_A)$ can be identified with  the symmetric group $S_n$ 
via
\begin{equation}\label{symmgr}S_n\ni w =  {\rm Diag}(r_{w(1)},\ldots,r_{w(n)}) =
(\H e_{w(1)},\ldots, \H e_{w(n)}),\end{equation} $ w\in S_n$
(see also equation (\ref{flnh})). All critical points are non-degenerate. 
\item[(ii)] Take $w\in S_n$ and $v:=s_{pq}w$,  for some $1\le p <q \le n$,
 where  
 $s_{pq}$ denotes the transposition of $p$ and $q$ in the symmetric group $S_n$. 
 Assume that  $h_A(w) >h_A(s_{pq}w)$.
 Then the subspace
 \begin{equation}\label{swpq}\S_{w,pq}:=K_{pq}.w,\end{equation}
 of $Fl_n(\H)$ is a metric sphere of dimension four in $({\mathcal H}_n, \langle \ , \ \rangle)$, for which
 $w$ and $s_{pq}w$ are the north, respectively south pole (with respect to the height function $h_A$). 
 Here $K_{pq}$ denotes the subgroup of 
$Sp(n)$ consisting of matrices with all entries zero, except for those on the diagonal and on the positions  $pq$ and  $qp$. The meridians of $\S_{w,pq}$  are gradient lines  between $w$ and $s_{pq}w$ for the function $h_A: Fl_n(\H)\to \bR$ with respect to the submanifold metric induced by $\langle \ , \ \rangle$.  
\item[(iii)] The negative space of the Hessian of $h_A$ at $w$ is
$\bigoplus_{(p,q)\in {\mathcal I}} T_w\S_{w,pq}$, where by ${\mathcal I}$ we denote  
the set of all pairs $(p,q)$ with $1\le p<q \le n$ such that $h_A(w)>h_A(ws_{pq})$.\end{itemize}
\end{proposition}
 The sphere  $\S_{w,pq}$ can also be described as  
 the set of all flags $(L_1,\ldots,L_n)$ with the property  that
 $L_{\nu}=\H e_{\nu}$, for all $\nu \notin \{p,q\}$ and $L_p$ and $L_q$ are arbitrary (orthogonal) lines in $\H e_{w(p)}\oplus\H e_{w(q)}$.  It is obvious that $\S_{w,pq}$ can be identified with
 the projective line $\H P^1 =\P(\H e_{w(p)}\oplus \H e_{w(q)})$, which is indeed a four-sphere.

\section{The GKM type description of $H^*_K(Fl_n(\H))$}

We consider the canonical embedding of $U(1)=\{a+bi  \ | \ a^2+b^2=1\}$ into
 $Sp(1)=\{a+bi+cj+dk \ | \ a^2+b^2+c^2+d^2=1\}$. This induces an embedding of $T:=U(1)^n$ into $K=Sp(1)^n$ (as spaces of diagonal matrices).  We are interested in the fixed points of the action of the groups $T$ and $K$ on $Fl_n(\H)$. 
 \begin{lemma}\label{lemmafirst} The fixed point sets $Fl_n(\H)^T$ and $Fl_n(\H)^K$ are both equal to $S_n$ (see the identification (\ref{symmgr})).  
\end{lemma}

\begin{proof}  It is obvious that  any flag of the form indicated in the lemma is fixed by
$K$. 
We  prove that  a flag fixed by $T$ has the form indicated in the theorem. Equivalently, we show that if  the vector $h=(h_1,\ldots,h_n) \in \H^n$ has the property that $L:=\H h$ is $T$-invariant, then $L$ must be $\H e_{\nu}$, for some $\nu\in\{1,\ldots,n\}$. Indeed, for any $A={\rm Diag}(z_1,\ldots,z_n)\in T$ we have $A.h \in L$, that is,  there exists $\lambda \in \H$ such that 
\begin{equation}\label{hbar}h_1\bar{z}_1=\lambda h_1,\ldots,  h_n\bar{z}_n=\lambda h_n.\end{equation} 
If $h_{\mu}\ne 0$ and $h_{\nu}\ne 0$ for $\nu \ne \mu$, we pick $z_{\mu}=1$ and
$z_{\nu}=-1$ and see that there exists no $\lambda$ satisfying (\ref{hbar}).
This finishes the proof. 
\end{proof}

We will  prove the following analogue of the Kirwan injectivity theorem.  

\begin{proposition}\label{propfirst} The map ${\imath}^*:H^*_K(Fl_n(\H)) \to H^*_K(S_n)$ induced by the inclusion $\imath : S_n \to Fl_n(\H)$ is injective.
\end{proposition} 

\begin{proof}  By Proposition \ref{morsecell}, $h_A$ is a Morse function, whose critical set is $S_n$. Let us order the latter set as 
$w_1,\ldots ,w_k$ such that $h_A(w_1)<h_A(w_2)< \ldots <h_A(w_k)$. 
Take $\epsilon >0$, smaller than the minimum  of $h_A(w_{\ell})-h_A(w_{\ell-1})$, where $2\le \ell \le k$.  
Denote $M_{\ell}:=h_A^{-1}(-\infty, h_a(w_{\ell})+\epsilon]$,
 and $S_n^{\ell} 
=S_n\cap M_{\ell}$. We prove by induction on 
$\ell\in \{1,\ldots, k\}$ that the map   
$\imath_{\ell}^* :H^*_K(M_{\ell}) \to H^*_K(S_n^{\ell})$ is injective. For $\ell=1$, the result is obviously true, as $M_{1}$ is equivariantly contractible to $\{w_1\}$. We assume the result is true for 
$\ell-1$. In the same way as in [To-We, diagram 2.5], we have the following commutative diagram. \begin{equation}\label{exactseq}
\vcenter{\xymatrix{
\cdots
\ar[r] &
H_{K}^*(M_{\ell},M_{\ell-1})
\ar[r]^{\textcircled{\small{2}}}\ar[d]^{\simeq} &
H_{K}^*(M_{\ell})
\ar[r] \ar[d]^{\textcircled{\small{1}}} &
H_{K}^*(M_{\ell-1})\ar[r] &  \cdots
\\
&
H_{K}^{*-{\rm index}(w_{\ell}) } (\{w_{\ell}\})
\ar[r]^{~~ \ \  \ \ \ \ \cup e_{\ell}}  &
H_{K}^*(\{w_{\ell}\}) &
\\
}}
\end{equation}  
Here $\simeq$ denotes the  isomorphism obtained by  composing the excision map $H^*_K(M_{\ell},M_{\ell-1})\simeq H^*_K(D,S)$ (where $D$, $S$ are the unit disk,
respectively unit sphere in the negative normal bundle of the Hessian of $h_A$ at the point $w_{\ell}$),  with the Thom isomorphism $H^*_K(D,S)\simeq 
 H^{*-{\rm index}(w_{\ell}) } _K(\{w_{\ell}\})$; the map $\textcircled{\small{1}}$ is induced by the inclusion of
 $\{w_{\ell}\}$ in $M_{\ell}$; the class 
  $e_{\ell}\in H_K^{{\rm index}(w_{\ell})}(\{w_{\ell}\})=H^{{\rm index}(w_{\ell})}(BK)$ is the $K$-equivariant Euler class of the negative space of the Hessian of $h_A$ at 
  $w_{\ell}$. 
According to Lemma \ref{lemmafirst}, $w_{\ell}$ is an isolated fixed point of the $T$ action on $Fl_n(\H)$. By the Atiyah-Bott lemma (see  appendix A of our paper) we deduce that the multiplication by 
$e_{\ell}$ is an injective map. We deduce that the long exact sequence of the pair 
$(M_{\ell},M_{\ell-1})$ splits into short exact sequences of the form
 $$
0\longrightarrow H_{K}^*(M_{\ell},M_{\ell-1})
\longrightarrow
H_{K}^*(M_{\ell})
\longrightarrow 
H_{K}^*(M_{\ell-1})\longrightarrow 0.$$
Let us consider now the following commutative diagram.
$$
\begin{array}{c}
\xymatrix{
0  \ar[r] & H_K^*(M_{\ell},M_{\ell-1}) \ar[r]\ar[d]_{\textcircled{\small{3}}}
& H_K^*(M_{\ell}) \ar[r]\ar[d]_{\imath_{\ell}^*} & H_K^*(M_{\ell-1}) 
\ar[r]\ar[d]_{\imath_{\ell-1}^*}
& 0 \\
0  \ar[r] & H_K^*(\{w_{\ell}\})\ar[r] &
H_K^*(S_n^{\ell}) \ar[r] & H_K^*(S_n^{\ell-1}) \ar[r]
& 0 }
\end{array}
$$
  where we have identified $H_K^*(S_n^{\ell},S_n^{\ell-1})=H^*_K(\{w_{\ell}\})$.   
 The map $\imath_{\ell-1}^*$ is injective, by the induction hypothesis. 
The map $\textcircled{\small{3}}$ is the same as the composition of $\textcircled{\small 1}$ and
$\textcircled{\small 2}$ (see diagram (\ref{exactseq})), thus it is injective as well. By a diagram chase we deduce that $\imath_{\ell}^*$ is injective.   
  \end{proof}  
 
  The following result will be used later.

\begin{lemma}\label{pqqp} The group $$G:=(\H^*)^n=\{{\rm Diag}(\gamma_1,\ldots,\gamma_n) \ | \ \gamma_{\nu}\in \H^*\}$$
acts transitively on $\S_{w,pq}\setminus \{w,ws_{pq}\}$ (see Proposition \ref{morsecell} for the definition of $\S_{w,pq}$). There exists a point in the latter space whose stabilizer is the group
$$G_{pq}=\{(\gamma_1,\ldots,\gamma_n)\in G \ | \ \gamma_p=\gamma_q\}.$$
\end{lemma}

\begin{proof} We use the description of $\S_{w,pq}$ given  after 
Proposition \ref{morsecell}. To prove the lemma, we only need to note  that for any two lines
$\H v_1$ and $\H v_2$  in $\H \oplus \H$  which are different from the ``coordinate
axes" $\H e_1$ and
$ \H e_2$ we have
$$\H v_1 = \H v_2 \cdot    \left(
\begin{array}{cc}
   \gamma_1 & 0 \\
0 &  \gamma_2  \\
\end{array}
\right)
$$
for some $\gamma_1,\gamma_2\in \H^*$. Moreover, in the same set-up, the 
$\H^*\times \H^*$ stabilizer of the line $\H(e_1+e_2)$ is the group
$\{(\gamma,\gamma) \ | \ \gamma\in \H^*\}.$

\end{proof}

Set  $$N:=\bigcup \S_{w,pq},$$
where the union runs over all $w,p,q$ such that $w\in S_n$, $1\le p <q \le n$, and $h_A(w)>h_A(ws_{pq})$. 
We  recall that since $K=Sp(1)^n$,  the $K$-equivariant cohomology of a point pt. is given by 
$$H^*_K({\rm pt.})=H^*(BK)=H^*(BSp(1)^n)=\Z[u_1,\ldots,u_n].$$

\begin{lemma}\label{separat} (i) Fix $p\in \{1,\ldots,n\}$ and consider the action of $K$  on  the vector bundle   $\H$ over a point {\rm pt.} given by $\gamma.h:=\gamma_{p}h,$
for any $\gamma=(\gamma_1,\ldots,\gamma_n)\in K$ and any $h\in \H$. 
The $K$-equivariant Euler class of the bundle is $u_{p}$. 

(ii) Fix $p, q\in \{1,\ldots,n\}$ and consider the action of $K$  on  the vector bundle   
$\H$ over a point {\rm pt.} given by $\gamma.h:=\gamma_{p}h\bar{\gamma}_{q},$
for any $\gamma=(\gamma_1,\ldots,\gamma_n)\in K$ and any $h\in \H$. 
The $K$-equivariant Euler class of the bundle is $u_{p}-u_{q}$. 
 
\end{lemma}

\begin{proof} 
(i) 
 Let $E=EK=E(Sp(1))^n$ be the total space of the classifying bundle of $K$.      
 By  definition (see for instance section 6.2 of [Gu-Gi-Ka]), the equivariant Euler class is  $$e_K(\H)=e(E \times_K \H),$$ where $E \times_K \H$ is a vector bundle with fiber 
 $\H$ over $E/K=BK=B(Sp(1))^n$. Thus
 $e_K( \H)$ is an element of degree 4 in $H^*(BK)=
 \bR[u_1,\ldots,u_n]$, hence a degree one polynomial in $u_1,\ldots,u_n$.    
We denote $BK=B(Sp(1))^n=\H P^{\infty}_1\times \ldots \H P^{\infty}_n$, where
$\H P^{\infty}_{\nu}$, $1\le \nu \le n$, are copies of $\H P^{\infty}$
(see e.g. [Hu, chapter 8, Theorem 6.1]). The latter is the space of all $\H$ lines
(that is, one dimensional 
$\H$-submodules) in $\H^{\infty}$, where the scalar multiplication is {\it from the right}. By this we mean that if $e\in \H^{\infty}$ and $h\in \H$, then, by definition, $h.e:=e\bar{h}$. 
   It is sufficient to show that if 
$\imath_{\nu}$ is the inclusion of $\H P^{\infty}_{\nu}$ into $BK$, then we have
$$\imath_{\nu}^*(e(E \times_K \H))=
\begin{cases} 0,\ {\rm  if} \  \nu\neq p\\
x, \ {\rm if} \  \nu=p.
\end{cases}
 $$
Here $x\in H^4(\H P^{\infty})$ denotes the Euler class of the tautological vector bundle    
  over $\H P^{\infty}$.  To show this, we note first that if $\nu\neq p$, then the restriction of the bundle $E\times_K \H$ to $\H P^{\infty}_\nu$ is the trivial bundle. Also, the restriction of 
  $E \times_K \H$ to $\H P^{\infty}_p$ is the bundle $ESp(1) \times_{Sp(1)} \H$ over 
  $ESp(1)/Sp(1)=\H P^{\infty}$ (here the action of $Sp(1)$ on $\H$ is by left multiplication). The latter bundle is actually the tautological bundle
   $$\tau=\{(v,L) \ | \ L \ {\rm is  \ a \ } \H \ {\rm  line \ in \ } \H^{\infty}, v\in L\}$$  
  over $\H P^{\infty}$.  Indeed, we take into account that $ESp(1)$ is the unit sphere $S^{\infty}$ in $\H^{\infty}$ (again by [Hu, chapter 8, Theorem 6.1]), thus we can identify $ESp(1)\times_{Sp(1)}\H$ with $\tau$ via 
  $$[e,h]\mapsto (eh, e\H ),$$
 for all $e\in ESp(1)$ and $h\in \H$. 
 We conclude that $\imath_p^*(e(E \times_K \H))=e(\tau)=x$.
  
  (ii) This time we have to check  that
 $$\imath_{\nu}^*(e(E \times_K \H))=
\begin{cases} 0,\ {\rm  if} \  \nu\neq p \ {\rm or} \ q\\
x, \ {\rm if} \  \nu=p\\
-x, \ {\rm if} \ \nu=q.
\end{cases}
 $$ 
 The cases $\nu\neq q$ have been discussed before, at point (i). 
 To analyze the case $\nu=q$, we note that the restriction of $E \times_K \H$ to $\H P^{\infty}_q$ is $ESp(1)\times_{Sp(1)}\H$, where $Sp(1)$ acts on $\H$ from the right. 
The map $$ESp(1)\times_{Sp(1)}\H\ni [e,h]\mapsto (e\bar{h},e\H)$$
is an isomorphism between $ESp(1)\times_{Sp(1)}\H$ and   the  rank four vector bundle over $\H P^{\infty}$,  call it $\bar{\tau}$, whose fiber over $L$ is $\bar{L}:=\{\bar{v} \ | \ v\in L\}.$
Because the linear automorphism of  $\bR^4 =\H$  given by
$$(x_1,x_2,x_3,x_4) =x_1+x_2i+x_3j+x_4k\mapsto  
x_1-x_2i-x_3j-x_4k = (x_1,-x_2,-x_3,-x_4)$$
is orientation changing (its determinant is equal to $-1$), we deduce 
that $$\imath_q^*(e(E \times_K \H))=e(\bar{\tau})=-e(\tau)=-x.$$
The claim, and also the lemma, are completely proved.
\end{proof} 

The following result will also be needed later.

\begin{lemma}\label{lemmaimp} Fix $w\in S_n$, $\epsilon >0$ strictly smaller than 
$|h_A(w)-h_A(v)|$, for any $v\in S_n$, $v\neq w$, and set
$N_+:=N\cap h_A^{-1}(-\infty,h_A(w)+\epsilon]$, $S_n^-=S_n\cap h_A^{-1}(-\infty,h_A(w)-\epsilon]$. Let $\eta$ be an element of $H^*_K(N^+)$ which vanishes when restricted to $S_n^-$. Then the restriction $\eta|_{w}$ is a multiple of $e_w$, where $e_w$ is the $K$-equivariant Euler class of the negative space of the Hessian of 
$h_A$ at $w$.
\end{lemma}  

\begin{proof} By Proposition \ref{morsecell}, (iii), the negative space of the Hessian of $h_A$ at $w$ is 
 $\bigoplus_{(p,q)\in {\mathcal I}} T_w\S_{w,pq}$.  For each $(p,q)\in {\mathcal I}$, the class $\eta_1:=\eta|_{\S_{w,pq}}\in H^*_K(\S_{w,pq})$  vanishes when restricted to $ws_{pq}$, which is the South pole of $\S_{w,pq}$. 
Consequently, $\eta_1$ vanishes when restricted to $h_A^{-1}(-\infty,h_A(w)-\epsilon]\cap \S_{w,pq}$,
since the latter can be equivariantly retracted onto $ws_{pq}$. The long exact sequence of the pair 
$(\S_{w,pq}, h_A^{-1}(-\infty,h_A(w)-\epsilon]\cap \S_{w,pq})$ is
$$\ldots \to H^*_K(\S_{w,pq}, h_A^{-1}(-\infty,h_A(w)-\epsilon]\cap \S_{w,pq})
\to H^*_K(\S_{w,pq})\to H^*_K(h_A^{-1}(-\infty,h_A(w)-\epsilon]\cap \S_{w,pq})\to \ldots. $$
Here we can replace $H^*_K(\S_{w,pq}, h_A^{-1}(-\infty,h_A(w)-\epsilon]\cap \S_{w,pq})$ with 
$H^*_K(D,S)$ where $D,S$ are a disk, respectively sphere, around $w$ (the north pole of 
$\S_{w,pq}$). In turn, $H^*_K(D,S)\simeq H^{*-4}_K(\{w\})$. Like in the proof of  
Proposition \ref{propfirst} (see diagram (\ref{exactseq})), we deduce that 
$\eta_1|_w$, which is the same as $\eta|_w$, is a multiple of $e_K(T_w(\S_{w,pq}))$,  the $K$-equivariant Euler class of 
the tangent space $T_w(\S_{w,pq})$. 
Now we recall that $w$ is a fixed point of the $K$ action and
 $\S_{w,pq}=K_{pq}.w$; consequently,  the tangent space at $w$ to 
$\S_{w,pq}$ is $$T_w(\S_{w,pq})=\k^0_{pq}.w$$
where the dot indicates the infinitesimal action and 
$$\k^0_{pq}=\{h E_{pq}+\bar{h}E_{qp} \ | \ h \in \H\}.$$ 
Here $E_{pq}$ denotes the $n\times n$ matrix whose entries are all 0, except for the one on position $pq$, which is equal to 1 (and the same for $E_{qp}$).  
This implies that $T_w(\S_{w,pq})$ is  $K$-equivariantly isomorphic to
 $\H$, with the $K$ action given by 
 $$(\gamma_1,\ldots,\gamma_n).h:=\gamma_ph{\gamma_q}^{-1},$$
 for any $(\gamma_1,\ldots,\gamma_n)\in K=(Sp(1))^n$ and any $h\in \H$. 
 By Lemma \ref{separat}, we have
 $$e_K(T_w(\S_{w,pq}))=u_p-u_q.$$
 
Thus the polynomial $\eta|_w$ is divisible by 
$u_p-u_q$, for all $(p,q)\in \I$. Because the latter polynomials are relatively prime with each other, we deduce that $\eta|_w$ is actually divisible by their product, which is just
$e_w$.   
\end{proof}

\begin{lemma}\label{same} If $\imath$, $\jmath$ are the inclusion maps of $S_n$ into $Fl_n(\H)$, respectively $N$, then the images of  $\imath^*:H_K^*(Fl_n(\H))\to H_K^*(S_n)$ and $\jmath^*:H_K^*(N)\to H_K^*(S_n)$ are the same.
\end{lemma}
 \begin{proof} Like in the proof of  Proposition \ref{propfirst}, we order $S_n=\{w_1,\ldots ,w_k\}$ such that  $h_A(w_1)<h_A(w_2)< \ldots <h_A(w_k)$. We choose $\epsilon >0$ smaller than the minimum of $h_A(w_{\ell})-h_A(w_{\ell-1})$, where $\ell=2,\ldots, k$. 
  We denote $M_{\ell}=h_A^{-1}(-\infty, h_A(w_{\ell})+\epsilon])$, $N_{\ell}:= N\cap M_{\ell}$,
  $S_n^{\ell}=S_n\cap M_{\ell}$  and will prove by induction on $\ell$ 
 that the images of the maps   
$\imath_{\ell}^*:H^*_K(M_{\ell}) \to H^*_K(S_n^{\ell})$ and $\jmath_{\ell}^*: H^*_K(N_{\ell}) \to H^*_K(S_n^{\ell})$ are the same.   For $\ell =1$, the assertion is clear, because both $M_1$ and its subset   $N_1$ can be retracted equivariantly onto $w_1$. Now we assume that the assertion is true for $\ell -1$ and we prove that it is true for $\ell$. Let us consider the following commutative diagram. 
\begin{equation}\label{dia1} \begin{array}{c}
\xymatrix{
 H_K^*(N_{{\ell}}) \ar[r]\ar[d]_{\jmath_{\ell }^*} & H_K^*(N_{\ell-1}) 
\ar[d]_{\jmath_{\ell -1}^*}
\\
H_K^*(S_n^{\ell })\ar[r]_{r}&
H_K^*(S_n^{\ell -1} ) }
\end{array}
\end{equation}
where $r$ is the restriction map.  We deduce that $r$ maps 
${\rm im}(\jmath_{\ell}^*)$ to ${\rm im}(\jmath_{\ell-1}^*)$.  From now on, by $r$ we will denote the induced map
 $$r:  {\rm im}(\jmath_{\ell}^*)\to {\rm im}(\jmath_{\ell-1}^*).$$
Let us consider another commutative diagram, namely
$$
\begin{array}{c}
\xymatrix{
0  \ar[r] & H_K^*(M_{\ell},M_{\ell-1}) \ar[r]\ar[d]_{h}
& H_K^*(M_{\ell}) \ar[r]\ar[d]_{\imath_{\ell}^*} & H_K^*(M_{\ell-1}) 
\ar[r]\ar[d]_{\imath_{\ell -1}^*}
& 0 \\
0  \ar[r] & \ker r\ar[r]_{g} &
{\rm im}(\jmath_{\ell}^*) \ar[r]_{r} & {\rm im}(\jmath_{\ell-1}^*)  \ar[r]
& 0 }
\end{array}
$$
In this diagram we have made two identifications: first, $H^*_K(S_n^{\ell},S_n^{\ell-1})$  is identified with $H^*_K(\{w_{\ell}\})$; then $H^*_K(\{w_{\ell}\})$ is canonically embedded in $H^*_K(S_n^{\ell})$. In this way,  $h$ goes from $H^*_K(M_{\ell},M_{\ell-1})$ to $H^*_K(S_n^{\ell},S_n^{\ell-1})=H^*_K(\{w_{\ell}\})$; but one can easily see from the diagram that $h$ actually takes  values in the subspace $\ker r$ of the latter space. The map $g$ is just the inclusion map. 
 We will prove that the image of $h$ is the whole $\ker r$. To this end, we take 
 $\eta\in H^*_K(N_{\ell})$ such that $r(\jmath^*_{\ell}(\eta))=0$ and we prove that $\jmath_{\ell}^*(\eta)$ is in the image of $g\circ h$, or, equivalently, that $\eta|_{w_{\ell}}$ is in the image of 
 $h$. 
 From the commutative diagram (\ref{dia1}) we deduce that the restriction of $\eta$ to $S_n^{\ell -1}$ is equal to 0. From Lemma \ref{lemmaimp} we deduce that 
 $\eta|_{w_{\ell}}$ is a multiple of the Euler class $e_{w_{\ell}}$.  
 Now let us consider again the diagram (\ref{exactseq}). We deduce that $\eta|_{w_{\ell}}$ is in the image of $\textcircled{\small{1}}\circ \textcircled{\small{2}}$, which is the same as $h$. In conclusion, we can use that $h$ and $\imath_{{\ell}-1}^*$ are surjective and, by a chase diagram,  deduce   that $\imath_{\ell}^*$ is surjective as well. The proposition is proved.
  \end{proof}
  
Now we are ready to characterize the image of $\imath^*$.   
  
\begin{proposition}\label{imagei} The image of $\imath^*:H^*_K(Fl_n(\H)) \to H^*_K(S_n)=
\bigoplus_{w\in S_n}\bR[u_1,\ldots,u_n]$ is
\begin{equation}\label{image} \{(f_w)_{w\in S_n} \ | \ u_p-u_q \ {\rm divides} \ \ f_w-f_{ws_{pq}}, \forall w\in S_n,\forall 1\le p <q \le n\}.\end{equation}
\end{proposition}

\begin{proof} Denote by $\mathfrak{im}$ the space described by (\ref{image}). By Lemma \ref{same}, it is sufficient to show that the image of $\jmath^*:H^*_K(N)\to H^*_K(S_n)$ is equal to $\mathfrak{im}$. 
First, let $(f_w)_{w\in S_n}$ be in the image of $\jmath^*$. Pick $w\in S_n$ and $p,q$ integers such that $1\le p <q \le n$ and $h_A(w)>h_A(ws_{pq})$. The pair $(f_w,f_{ws_{pq}})$ is in the image of the restriction map $H^*_K(\S_{w,pq})\to H^*_K(\{w,ws_{pq}\}).$ 

\noindent {\it Claim.} The polynomial $f_w-f_{ws_{pq}}$ is divisible by $u_p-u_q$.

 Put $U_1:=\S_{w,pq}\setminus \{w\}$, $U_2:=\S_{w,pq}\setminus \{ws_{pq}\}$. We use the Mayer-Vietoris sequence in equivariant cohomology for the pair $(U_1,U_2)$. We note that $U_1\cap U_2=\S_{w,pq}\setminus \{w,ws_{pq}\}$.
 By Lemma \ref{pqqp}, we can identify $\S_{w,pq}\setminus \{w,ws_{pq}\}$ with the homogeneous space $G/G_{pq}$. 
Now, because $G/K$ is  contractible, we have  that $$H^*_K(U_1\cap U_2) = H^*_{G}(G/G_{pq}).$$   The ring in the right hand side consists of all polynomials  in $u_1,\ldots,u_n$ which are divisible by $u_p-u_q$, that is 
$$H^*_K(U_1\cap U_2)=\Z[u_1,\ldots,u_n]/\langle u_p-u_q\rangle.$$ 
Also note that since $U_1$ and $U_2$ are contractible, we have $H^*_K(U_1)=H^*_K(\{w\})$ and 
$H^*_K(U_2)=H^*_K(\{ws_{pq}\})$, which are both isomorphic to 
$H^*(BK)=\bR[u_1,\ldots,u_n]$. Moreover, the restriction maps 
$H^*_K(U_1)\to H^*_K(U_1\cap U_2)$ and $H^*_K(U_2)\to H^*_K(U_1\cap U_2)$ are both the canonical projection 
$$\Z[u_1,\ldots,u_n]\to \Z[u_1,\ldots,u_n]/\langle u_p-u_q\rangle.$$
The Mayer-Vietoris sequence is
\begin{equation}\label{mv}\ldots \to H_K^*(\S_{w,pq}) \to H_K^*(U_1)\oplus H_K^*(U_2) \stackrel{g}{\to} H_K^*(U_1\cap U_2) \to \ldots .\end{equation}
The claim follows from the fact that $g(f_{w},f_{ws_{pq}})=  f_{w}-f_{ws_{pq}}$ mod $\langle u_p-u_q\rangle$ is equal to 0.

Now we prove that $\mathfrak{im}$ is contained in the image of $\jmath^*$.  Let
$(f_w)_{w\in S_n}$ be an element of $\mathfrak{im}$. 
From the exact sequence $(\ref{mv})$, we deduce that for each $w\in S_n$ and each pair $p,q$ with $1\le p <q \le n$, $h_A(w)>h_A(ws_{pq})$, there exists $\alpha_{w,pq}\in H^*_K(\S_{w,pq})$ with $\alpha_{w,pq}|_w=f_w$ and $\alpha_{w,pq}|_{ws_{pq}} = f_{ws_{pq}}$. 
A simple argument (using again a Mayer-Vietoris sequence) shows that there exists 
$\alpha \in H^*_K(N)$ such that 
$\alpha|_{\S_{w,pq}}=\alpha_{w,pq}$. This implies $(f_w)_{w\in S_n}=\jmath^*(\alpha)$.
\end{proof} 
 
 \noindent {\bf Remarks.}  1. It is likely that the main results of this section, namely
 Proposition \ref{propfirst} and Proposition \ref{imagei}, can be proved with the methods of [Ha-He-Ho]. 
 
 2. A GKM type description similar to the one given in Proposition \ref{imagei}, for
 $(\Z/ 2\Z)^n$-equivariant cohomology with coefficients in $\Z/2\Z$ can be deduced 
  from [Bi-Ho-Gu,  Theorem C] (we recall that $Fl_n(\H)$ can be regarded as the real locus of a coadjoint orbit of the group $SU(2n)$). 
  
The following result will complete the proof of Theorem \ref{main}.
 We consider the tautological vector bundles 
${\mathcal L}_i$ over $Fl_n(\H)$, $1\le i\le n$. That is, the fibre of 
${\mathcal L}_i$ over $(L_1,\ldots,L_n)\in Fl_n(\H)$ is $L_i$. 

\begin{lemma}\label{uwi} Take $w\in S_n=(Fl_n(\H)^K)$, and identify
 $H^*_K(\{w\})=H^*(BK)=\Z[u_1,\ldots,u_n]$.  
Then the equivariant Euler class $e_K({\mathcal L}_{\nu})$ restricted to $w$ is equal to 
$u_{w(\nu)}$, for all $\nu\in\{1,\ldots,n\}$.
\end{lemma}

\begin{proof} The cohomology class  $e_K({\mathcal L}_{\nu})|_{w}$ is the equivariant  Euler class of the space ${\mathcal L}_{\nu}|_{w}=\H $ with the $K$ action 
$$(\gamma_1,\ldots,\gamma_n).h=\gamma_{w(\nu)}h,$$
for all $(\gamma_1,\ldots,\gamma_n)\in K= Sp(1)^n$ and $h\in \H$. The result follows from 
Lemma \ref{separat}. \end{proof}  

Now we are ready to prove our main result.

\noindent {\it Proof of Theorem \ref{main}.} Let $Fl_n(\bC)$ be the space of flags in $\bC^n$, which can be equipped in a natural way with the action of the torus
$T:=(S^1)^n$. The idea of the proof is to compare the equivariant cohomology rings $H^*_T(Fl_n(\bC))$ and $H^*_K(Fl_n(\H))$. The first one can be computed using the GKM theory (cf. [Ho-Gu-Za]), as follows.
 Like in the quaternionic case explained here,  the fixed point set  $Fl_n(\bC)^T$ can be identified with the symmetric group $S_n$, and the restriction map $H^*_T(Fl_n(\bC))\to H^*_T(S_n)$ is injective. Moreover, the image of the latter map consists of all sequences of polynomials
  $(f_w)_{w\in S_n}$ such that
 $f_w-f_{s_{pq}w}$ is divisible by $\tilde{u}_p-\tilde{u}_q$, for all $w\in S_n$ and $1\le p<q\le n$,
 where we have identified $$H^*(BT)=\Z[\tilde{u}_1,\ldots,\tilde{u}_n]$$ (compare to Proposition \ref{imagei}). 
 The equivariant Euler classes $e_T(\tilde{\mathcal L}_{\nu})$ of the tautological
 complex line  bundles over $Fl_n(\bC)$ have the property that 
 $e_T(\tilde{\mathcal L}_{\nu})|_w=\tilde{u}_{w(\nu)}$, for all $1\le \nu \le n$ (compare to Lemma \ref{uwi}).
 On the other hand, we have the Borel-type description of $H^*_T(Fl_n(\bC))$, namely
  $$H^*_T(Fl_n(\bC)) \simeq \Z[x_1,\ldots,x_n,\tilde{u}_1,\ldots,\tilde{u}_n]
/\langle (1+x_1)\ldots (1+x_n)=(1+\tilde{u}_1)\ldots (1+\tilde{u}_n)\rangle,$$
via $e_T(\tilde{\mathcal L}_{\nu}) \mapsto x_{\nu}$, $1\le \nu \le n$.  Theorem \ref{main}  follows.
$\hfill \square$  
  
\section{Appendix A. The Atiyah-Bott lemma}
We will prove the following version  of the
Atiyah and Bott's [At-Bo] lemma. We recall that all cohomology rings are with integer
coefficients (unless otherwise specified). 
 \begin{lemma}\label{lemmalast}  Let $V$ be an even-dimensional real vector space,
 with the linear action of the group $K:=Sp(1)^n$.  Assume that the only fixed point of the action of $T:=(S^1)^n\subset (Sp(1))^n$ on $V$ is $0$.
 If one regards $V$ as a vector bundle over the point $0$, then
 the equivariant Euler class $e_K(V) \in H^*_K(\{0\})=H^*(BK)$ is different from zero
 (hence it is not a zero divisor).
 \end{lemma}
 \begin{proof} 
 The inclusion $H^*_K(\{0\})\to H^*_T(\{0\})$ maps 
 $e_K(V)$ to $e_T(V)$ (because if $E:=EK=ET$, then the line bundle
 $E\times_T V$ over $E/T$ is the pullback of $E\times_K V$ over $E/K$ via the natural map $E/T \to E/K$). 
 On the other hand, the natural map $H^*_K(\{0\})\to H^*_T(\{0\})$ is injective.
 The reasons are as follows: first,  the map $H^*(BK,\Q)=H^*_K(\{0\},\Q)\to H^*_T(\{0\},\Q)
 =H^*(BT,\Q)$ is injective (its image is actually the set of all elements invariant under the Weyl group action); second,  both $H^*(BK)=
 H^*((\H P^{\infty})^n)$ and $H^*(BT)=H^*((\bC P^{\infty})^n)$ are torsion free.
 So it is sufficient to show that $e_T(V)\in H^*_T(\{0\})=\bR[u_1,\ldots,u_n]$  is different from zero.   Since the representation of $T$ on $V$ has no nonzero fixed points, we have $V=\oplus_{i=1}^mL_i$, where $L_i$ are $1$-dimensional  complex representations of $T$. Thus we have
$$e_{T}(V)=c_m^{T}(V)=c_m^{T}(\oplus_{i=1}^mL_i)
=c_1^{T}(L_1)\ldots c_1^{T}(L_m),$$
where $c_m^{T}$ and $c_1^{T}$ denote the $T$-equivariant Chern classes.
Each Chern class $c_1^{T}(L_i)$ is different from zero, since the 1-dimensional complex  representations of $T$ are labeled  by the character group 
${\rm Hom}(T,S^1)$, and the map ${\rm Hom}(T,S^1)\to H^2(BT)$ given by
$L\mapsto c_1^{T}(L)$ is a linear isomorphism (see for instance [Hu, chapter 20, section 11]). This finishes the proof.
\end{proof}

\section{Appendix B. Height functions on isoparametric submanifolds}
  
The goal of this section is to provide a proof of Proposition \ref{morsecell}, and also to achieve a better understanding of the spheres $\S_{w,pq}$, which are important objects of our paper.  We will place ourselves in the more general context of isoparametric submanifolds.
We recall (see for instance [Pa-Te, chapter 6]) that   an $n$-dimensional submanifold   $M\subset \bR^{n+k}$  which is  closed, complete with respect to the induced metric,
and full
(i.e. not contained in any affine subspace) is called 
  {\it isoparametric}  if  any normal vector at a point
of $M$
 can be extended
 to a parallel normal vector field $\xi$ on $M$ with the property that the eigenvalues of the shape operators
$A_{\xi(x)}$ (i.e. the principal curvatures) are independent of $x\in M$, as values and multiplicities.
It follows that for  $x\in M$, the set  $\{A_{\xi(x)}: \xi(x) \in \nu (M)_x\}$ is a commutative family of
selfadjoint endomorphisms of $T_x(M)$, and so it determines a decomposition of 
$T_x(M)$  as a
direct sum of common eigenspaces
$E_1(x), E_2(x), ...,E_r(x)$. There exist normal vectors  $\eta_1(x), \eta_2(x),..., \eta_r(x)$
such that
$$A_{\xi(x)}|_{E_i(x)}=\langle \xi(x), \eta_i(x) \rangle {\rm id}_{E_i(x)},$$
for all $\xi(x) \in \nu_x (M)$, $1\le i\le r$. By parallel extension in the normal bundle we obtain
the vector fields $\eta_1, \ldots, \eta_r$.
The eigenspaces from above give rise to the distributions $E_1, \ldots, E_r$ on $M$, which are called
the  {\it  curvature distributions}.
The numbers $$m_i= {\rm rank} E_i,$$ $1\leq i \leq r$, are  the {\em
multiplicities } of $M$. We fix a point $x_0\in M$ and we consider the normal space
$$\nu_0:=\nu_{x_0}(M).$$
In the affine space $x_0+\nu_0$ we  consider the hyperplanes  
$$\ell_i(x_0):= \{x_0+\xi(x_0) : \xi(x_0)\in \nu_0, \langle \eta_i(x_0), \xi(x_0) \rangle =1\},$$
$1\le i \le r$; one can show that they have a unique intersection point, call
it $c_0$, which is independent of the choice of $x_0$. Moreover, $M$ is contained
in a sphere with center at $c_0$. 
We do not lose any generality if we
assume that $M$ is contained in the unit sphere $S^{n+k-1}$, hence
$c_0$ is just the origin $0$ and $x_0+\nu_0=\nu_0$ (because $x_0\in \nu_0$).  One shows that the group of linear
transformations of $\nu_0$ generated by the reflections about
$\ell_1(x_0),\ldots, \ell_r(x_0)$ is a Coxeter group. We denote it by $W$ and call it the
{\it Weyl group} of $M$. We have $$\nu_0\cap M = W.x_0.$$

For each $i\in\{1,\ldots,r\}$, the distribution $E_i$ is integrable and the leaf through $x\in M$ of $E_i$ is a round distance  sphere $S_i(x)$, of dimension $m_i$, whose center is the orthogonal projection of $x$ on $\ell_i(x)$. These are called the {\it curvature spheres}.
We note that $S_i(x)$ contains $x$ and $s_ix$ as antipodal points.   We have the following result.

\begin{proposition}\label{isopara} Let $a \in \nu_0$  contained in (the interior of) the same Weyl chamber as $x_0$ and let $h_a: M\to \bR$, $h_a(x)=\langle a,x\rangle$
be the corresponding height function. The following is true.
\begin{itemize}
\item[(a)] The gradient of $h_a$ at any $x\in M$ is  $a^{T_x(M)}$, that is, the orthogonal projection of $a$ on the tangent space $T_xM$.
\item[(b)] ${\rm Crit}(h_a)=W.x_0$ and the negative space of the hessian of $h_a$ at 
$x\in {\rm Crit}(h_a)$ is $\bigoplus_i E_i(x)$,  where $i$ runs over those indices in $\{1,\ldots, r\}$ with the property
  $h_a(x)>h_a(s_ix)$. 

\item[(c)] For any $x\in {\rm Crit} (h_a)$, and $i\in \{1,\ldots, n\}$ such that $h_a(x)>h_a(s_ix)$,  the meridians of $S_i(x)$ going from $x$ to $s_ix$   are gradient lines of the function $h_a: M \to \bR$ with respect to the submanifold metric on $M$.
  
\end{itemize}
\end{proposition}
\begin{proof} 
The points (a) and (b) are proved in [Pa-Te, chapter 6].
We only need to prove (c).  
This follows immediately from the fact that for any $z\in S_i(x)$, the vector $\nabla (h_a)(z)=a^{T_zM}$ is tangent
to $S_i(z)$ (because $a\in \nu_{x_0}(M) =\nu_x(M)$ is perpendicular to
$\bigoplus_{j\neq i}E_j(x) =\bigoplus_{j\neq i}E_j(z)$, cf.  [Pa-Te, Theorem 6.2.9
(iv) and  Proposition 6.2.6]). 
\end{proof}

Finally, Proposition \ref{morsecell} can be deduced from Proposition \ref{isopara} by taking into account that $Fl_n(\H)$ is an isoparametric submanifold of ${\mathcal H}_n$ with the 
following properties (cf. [Pa-Te, Example 6.5.6] for the symmetric space 
$SU(2n)/Sp(n)$):
\begin{itemize}
\item the normal space to $x_0:={\rm Diag}(r_1,\ldots,r_n)$ is $\nu_0=\d$
\item the Weyl group $W$ is $S_n$ acting in the obvious way on $\d$
\item  the curvature spheres through $wx_0$ are the orbits $K_{pq}.w$, where
$1\le p<q\le n$; thus all multiplicities are equal to 4.
\end{itemize}

\bibliographystyle{abbrv}

\end{document}